\documentclass[12pt,a4paper]{article}
\usepackage{amsfonts}
\usepackage{amssymb}
\usepackage{wasysym}
\usepackage{amscd}
\begin{document}

  \newtheorem{theor}{Theorem}[section]
  \newtheorem{prop}[theor]{Proposition}
  \newtheorem{cor}[theor]{Corollary}
  \newtheorem{lemma}[theor]{Lemma}
  \newtheorem{sublem}[theor]{Sublemma}
  \newtheorem{defin}[theor]{Definition}
  \newtheorem{conj}[theor]{Conjecture}
  \hfuzz2cm
  \def\deg{{\widehat {\rm deg}\,}}
  \def\bz{\mbox{\boldmath$\zeta$\unboldmath}}
  \def\bzs{\mbox{\boldmath$\zeta'$\unboldmath}}
  \def\bzo{\overline{\bz}}
  \def\odd{{\rm odd}}
  \def\a{\alpha}
  \def\bC{{\bf C}}
  \def\Td{{\rm Td}}
  \def\ch{{\rm ch}}
  \def\Hom{{\rm Hom}}
  \def\Ad{{\rm Ad}^{1,0}_{G/P}}
  \def\KO{H}
  \def\K1{K}
  \def\covol{{\rm covol}\,}
  \gdef\beginProof{\par{\bf Proof: }}
  \gdef\endProof{${\bf Q.E.D.}$\par}
  \def\a{\alpha}
  \def\b{\beta}
  \def\th{\Theta^E}
  \def\Tr{{\rm Tr}\,}
  \def\Alt{{\rm Alt}}
  \def\BN{{\bf N}}
  \def\BZ{{\bf Z}}
  \def\BQ{{\bf Q}}
  \def\BR{{\bf R}}
  \def\BC{{\bf C}}
  \def\BG{{\bf G}}
  \def\BH{{\bf H}}
  \def\SB{{\raise.1em\hbox to 0pt{\hskip.1em\slash\hss} S}}
  \def\Sp{{\bf Sp}}
  \def\GL{{\bf GL}}
  \def\O{{\bf O}}
  \def\S{\mathrm{Sym}\,}
  \def\L{\Lambda\,}
  \def\F{\mathcal{W}}
  \def\g{\mathfrak{g}}
  \def\k{\mathfrak{k}}
  \def\p{\mathfrak{p}}
  \def\i{\lrcorner}
  \def\tfrac#1#2{{\textstyle{#1\over#2}}}
\def\spentagon{{\pentagon}}
\def\lpentagon{{\large\pentagon}}
  \def\fiverm{}
  \author{Kai K\"ohler$^1$\\
Gregor Weingart$^2$}
  \title{Quaternionic analytic torsion}
\maketitle
  \begin{abstract}
   We define an (equivariant) quaternionic analytic torsion for antiselfdual
   vector bundles on quaternionic K\"ahler manifolds, using ideas by Leung
   and Yi. We compute this torsion for vector bundles on quaternionic
homogeneous spaces
   with respect to any isometry in the component of the identity, in terms of
roots and Weyl groups.
  \end{abstract}
  \begin{center}
  \end{center}
\footnotetext[1]{Centre de Math\'ematiques de Jussieu/C.P. 7012/2, place
Jussieu/F-75251 Paris Cedex 05/France/e-mail : koehler@math.jussieu.fr/URL:
http://www.math.jussieu.fr/$\tilde{\ }$koehler}
\footnotetext[2]{Mathematisches Institut/ Wegelerstr. 10/D-53115
Bonn/Germany/e-mail: gw@rhein.iam.uni-bonn.de/URL:
http://www.math.uni-bonn.de/people/gw}
\thispagestyle{empty}

\setcounter{page}{1}
\tableofcontents
\newpage
\parindent=0pt
\parskip=5pt

\section{Introduction}
  Analytic torsions were introduced by Ray and Singer as real numbers
  constructed using certain $\textbf{Z}$--graduated complexes of elliptic
  differential operators acting on forms with coefficients in vector
  bundles on compact manifolds. The real analytic torsion \cite{RS1} was defined
  for the de Rham-operator associated to flat Hermitian vector bundles
  on Riemannian manifolds. It was proven by Cheeger and M\"uller to equal
  a topological invariant, the Reidemeister torsion, which can be defined
  using a finite triangulation of the manifold. This implies that the real
  analytic torsion is a homeomorphy invariant which is not invariant under
  homotopy. Lott and Rothenberg pointed out that an equivariant version of
  this torsion still is a diffeomorphy invariant.

  The complex Ray-Singer (or holomorphic) torsion \cite{RS2} was defined for the
  Dolbeault--operator acting on antiholomorphic differential forms with
  coefficients in a holomorphic Hermitian vector bundle on a compact
  complex manifold. It turned out to play an important role in the
  Arakelov geometry of schemes over Dedekind rings. In fact it was
  shown by Bismut, Gillet and Soul\'e to provide a direct image in
  a K-theory of Hermitian vector bundles. This direct image verifies
  a Grothendieck-Riemann-Roch relation with Arakelov-Gillet-Soul\'e
  intersection theory, as was proven by Bismut, Lebeau, Gillet and
  Soul\'e. Later, K\"ohler and Roessler showed that an equivariant
  version of this direct image localizes on fixed point subschemes
  in Arakelov geometry. This had many applications in arithmetic
  geometry, algebra and global analysis.

  Thus it seems natural to investigate torsions for other
  $\textbf{Z}$--graded complexes occurring in geometry, in particular
  for quaternionic manifolds. A first step towards a definition
  of analytic torsion for quaternionic manifolds was made
  in an e--print by Leung and Yi \cite{ly}, using a complex first
  discussed by Salamon. We had problems understanding this very general, short and
ambiguous definition.
  In the present paper, we first give a thorough definition of an
  (equivariant) quaternionic torsion for quaternionic K\"ahler manifolds $M$,
  with coefficients in the antiselfdual vector bundles $\F$. This is done by
  carefully decomposing the action of a natural Dirac operator on Salamon's
  complex on these manifolds, i.e. on the complex
  $$
    0\quad\longrightarrow\quad\S^{k}H\otimes\F
     \quad\stackrel{d}\longrightarrow\quad\S^{k+1}H\otimes\L^{1,0}E^*\otimes\F
    \qquad\qquad
  $$
  \unskip
  $$
    \qquad\qquad\qquad\stackrel{d}\longrightarrow\quad\ldots\quad\stackrel{d}
    \longrightarrow\quad\S^{2n+k}H\!\otimes\L^{2n,0}E^*\otimes \F\quad
    \longrightarrow\quad0
  $$
  for a parameter $k\in\BN_0$ even and $TM\otimes_\BR\BC\cong H\otimes E$.
  The Laplace operator defining the torsion is the square of this Dirac
  operator. We detail the many traps to avoid in this construction.

  In the third section, we compute the equivariant quaternionic torsion for
  all known quaternionic K\"ahler manifolds of positive curvature, i.e.
  for the quaternionic homogeneous spaces of the compact type, with respect
  to the action of any element of the associated Lie group and any equivariant
  antiselfdual vector bundle. These spaces are known to be symmetric, and
  for any simple
  compact Lie group there is exactly one quaternionic homogeneous space.
  This computation proceeds very similar to previous computations of the
  real analytic torsion and the holomorphic torsion for all appropriate
  symmetric spaces by one of us. We regard this as further indication
  that the definition given here is a ``good'' definition of quaternionic
  torsion.

  For the real analytic torsion, this computation led to a homeomorphy
  classification of quotients of some odd-dimensional symmetric spaces
  of the compact type. For the holomorphic torsion, this computation gave
  evidence for the fixed point formula mentioned above. In combination
  with the fixed point formula in Arakelov geometry, it provided a new
  proof of the Jantzen sum formula classifying the lattice representations
  of Chevalley group schemes except for the cases $G_2, F_4, E_8$. Thus
  one can reasonably hope for interesting
  applications of our result. Remarkably, the formula for quaternionic
  torsion happens to have the very same structure as the formula for the
  holomorphic torsion on Hermitian symmetric spaces (thus, on different
  manifolds). In a forthcoming paper, we intend to relate the torsion to
  holomorphic torsion on the twistor space, which should as an application of the
computation done here
  provide a full proof of the Jantzen sum formula including the three
  exceptional cases.

  In the last section, we comment briefly on the special case of
  hyperk\"ahler manifolds, in which the quaternionic torsion can
  be expressed in terms of a Dolbeault--operator. Related work for this case has
  been done recently by Gerasimov and Kotov \cite{GK1}, \cite{GK2}.

  {\bf Acknowledgements: } We are indebted to Sebastian Goette for
  many fruitful discussions. The first author thanks the Deutsche
  Forschungsgemeinschaft for supporting him with a Heisenberg
  fellowship.

\section{Quaternionic analytic torsion}

  Perhaps the most fundamental difference between quaternionic geometry
  and complex geometry is the lack of a plausible notion of quaternionic
  differentiability, any such notion leads inevitably to a finite dimensional
  space of quaternionic differentiable functions even on $\BH^n$. Hence the
  stock of local transition functions is rather limited and it seems impossible
  to define a quaternionic manifold in terms of an atlas of holomorphic
  coordinate charts. Another way to express this difference between complex
  and quaternionic geometry is that each quaternionic manifold comes along
  with a distinguished ``projective'' equivalence class of torsion free
  connections, a feature unheard of in complex geometry but rather
  characteristic for so called parabolic geometries. In fact quaternionic
  geometry can be seen as an example for parabolic geometries and many
  of the aspects discussed below are more or less directly linked to
  this fact. The interested reader is referred to \cite{barker} for
  this point of view.

  A quaternionic manifold $M$ is a manifold of dimension $4n,\,n\geq 2,$
  endowed with a smooth quaternionic structure on its tangent spaces
  admitting an adapted torsion free connection. In other words $M$ is
  endowed with a reduction $\Sp(1)\cdot\GL_\BH(M)\subset\GL(M)$ of its
  frame bundle to the bundle of quaternionic frames with structure group
  $\Sp(1)\cdot\GL_n(\BH)\subset\GL_{4n}(\BR)$ tangent to some torsion
  free connection. The projective equivalence class of this connection
  is uniquely determined by the quaternionic structure in the sense that
  the adapted connections are parametrized by $1$--forms on $M$. A guiding
  principle in the construction of differential complexes on quaternionic
  manifolds is hence to twist with a trivialisable line bundle in order
  to make the differential operators independent of the choice of connection
  following Fegan's approach to the construction of conformally invariant
  differential operators \cite{feg}.

  Note that dimension 4 is explicitly excluded from the definition given
  above, in fact the group $\Sp(1)\cdot\GL_n(\BH),\,n=1,$ is exactly the
  conformal group and the existence of a torsion free connection imposes
  no integrability assumption whatsoever on the conformal structure.
  Consequently differential sequences like (\ref{dcx}) and (\ref{dscx})
  below fail in general to be complexes in conformal geometry. However
  there is a geometry in dimension 4 analogous to quaternionic geometry
  in higher dimensions $4n,\,n>1$, namely the so called half conformally
  flat geometry of conformal manifolds with vanishing self--dual Weyl tensor.
  Mutatis mutandis our considerations below are valid in half conformally
  flat geometry in dimension 4, in particular the differential sequences
  (\ref{dcx}) and (\ref{dscx}) become complexes under this integrability
  assumption.

  Any representation of the group $\Sp(1)\cdot\GL_n(\BH)$ gives rise to
  a vector bundle on $M$ associated to the quaternionic frame bundle
  $\Sp(1)\cdot\GL_\BH(M)$. Consider the two defining representations
  $\pi_{H}=\BC^2$ of $\Sp(1)$ and $\pi_{E}=\BC^{2n}$ of $\GL_n(\BH)$
  respectively, which both carry invariant quaternionic structures $J$
  by definition. Moreover the representation $\pi_{H}$ carries an invariant
  symplectic form $\sigma$, which is real (i.e. $\sigma(Jh_1,\,Jh_2)=
  \overline{\sigma(h_1,\,h_2)}$) and positive (i.e. $\sigma(h,\,Jh)>0$
  for all $h\neq 0$). The existence of a non--degenerate bilinear form
  $\sigma$ implies in particular that $\pi_H$ is equivalent to its dual
  $\pi_H^*$ as an $\Sp(1)$--representation via the musical isomorphism
  $\sharp:\;\pi_H\longrightarrow\pi_H^*,\,h\longmapsto h^\sharp$ or its
  inverse $\flat$ with $h^\sharp:=\sigma(h,\cdot)$.

  Notice that the complex determinant of an element of $\GL_n(\BH)\subset
  \GL_{2n}(\BC)$ is always a real positive number. Hence the representation
  $(\mathrm{det}\,\pi_E)^s$ of $\GL_n(\BH)$ is defined for any $s\in\BR$.
  Moreover $\mathbf{PGL}_n(\BH)\,:=\,\GL_n(\BH)/\BR^*$ is a real form of
  $\mathbf{SL}_{2n}(\BC)$ and thus all irreducible representations
  of $\Sp(1)\times\GL_n(\BH)$ occur in tensor products of $\pi_H,\,\pi_E,\,
  \pi_E^*$ and $(\mathrm{det}\,\pi_E)^s$ with $s\in\BR$. The irreducible
  representations occuring in a tensor product $(\mathrm{det}\,\pi_E)^s
  \otimes\pi_{H}^{\otimes k}\otimes\pi_{E}^{\otimes a}\otimes\pi_{E}^{*\otimes
  b}$ with $k+a+b$ even descend to $\Sp(1)\cdot\GL_n(\BH)$, in particular
  all irreducible representations of $\Sp(1)\cdot\GL_n(\BH)$ carry real
  structures and so do all vector bundles associated to the bundle of
  quaternionic frames e.~g.~the complexified tangent bundle $TM\otimes_\BR\BC$
  is associated to the representation $\pi_{H}\otimes\pi_{E}$. We will write
  $TM\otimes_\BR\BC\cong H\otimes E$ although this notation has to be
  taken with care as neither $E$ nor $H$ are globally defined vector
  bundles in general. The trivializable line bundles associated to the
  representations $(\mathrm{det}\,\pi_E)^s$ will be denoted by $L^s$.

  The invariant symplectic form on $\pi_{H}$ defines a real, positive section
  $\sigma_H$ of the vector bundle $\L^2H$, which is parallel for every adapted
  torsion free connection. Choosing similarly a real, positive section
  $\sigma_E$ of $\L^2E^*$ amounts to choosing a Riemannian metric on $M$
  compatible with the quaternionic structure. In particular quaternionic
  K\"ahler manifolds are quaternionic manifolds $M$ with a fixed real,
  positive section $\sigma_E$ of $\L^2E$, which is parallel for an adapted
  torsion free connection, necessarily equal to the Levi--Civita connection
  of the Riemannian metric $\sigma_H\otimes\sigma_E$.

  The tensor product decomposition $TM\otimes_\BR\BC\cong H\otimes E$
  of the complexified tangent bundle of a quaternionic manifold induces
  a corresponding decomposition $T^*M\otimes_\BR\BC\cong H\otimes E^*$
  of its cotangent bundle and of the whole exterior algebra of forms.
  According to the theory of Schur functors \cite{fh} this
  decomposition reads
  \begin{equation}\label{schur}
   \L^\bullet(T^*M\otimes_\BR\BC)\quad\cong\quad
   \bigoplus_{2n\geq a\geq b\geq 0\atop a+b\,=\,\bullet}
   \;\;\S^{a-b}H\,\otimes\,\L^{a,\,b}E^*\,,
  \end{equation}
  where $\L^{a,\,b}E^*\,\subset\,\L^aE^*\otimes\L^bE^*$ is the kernel of a
  $\GL_n(\BH)$--equivariant map
  $$
   0\;\;\longrightarrow\;\;\L^{a,\,b}E^*\;\;\stackrel\subset\longrightarrow
    \;\;\L^aE^*\,\otimes\,\L^bE^*\;\;\stackrel{\textrm{Pl}}\longrightarrow
    \;\;\L^{a+1}E^*\,\otimes\,\L^{b-1}E^*\,\;\;\longrightarrow\;\;0\,,
  $$
  whose precise definition is immaterial for the arguments below.
  Wedging with a $1$--form in $H\otimes E^*$ maps $\S^{a-b}H
  \otimes\L^{a,b}E^*$ to the sum $\S^{a-b+1}H\otimes\L^{a+1,b}E^*
  \oplus\S^{a-b-1}H\otimes\L^{a,b+1}E^*$ of course and a little more
  elaboration provides us with an explicit isomorphism (\ref{schur})
  such that
  \begin{eqnarray}
   (h\otimes\eta)\wedge &=&
   \frac1{a-b+1}\,h\cdot\otimes\eta\wedge\otimes\mathrm{id}
   \;+\;(-1)^{a-b}h^\sharp\i\otimes\mathrm{id}\otimes\eta\wedge\label{wf}\\
   && \qquad-\;\frac{(-1)^{a-b}}{a-b+1}(\mathrm{id}\otimes\textrm{Pl}^*)
      \circ(h^\sharp\i\otimes\eta\wedge\otimes\mathrm{id})\nonumber
  \end{eqnarray}
  with some linear map $\textrm{Pl}^*$ twin to $\textrm{Pl}$ above. Now
  for a quaternionic manifold the decomposition of the exterior algebra
  is respected by some torsion free connection and consequently the de
  Rham complex of $M$ gives rise both to a quotient complex and a
  subcomplex of the form:
  $$
   \begin{CD}
    0\!\rightarrow\;\;\quad\BC\;\;\quad\!\!@>d>>\!\!H\!\otimes\!\L^{1,0}E^*
    \!\!@>d>>\!\!\!\ldots \!\!\!@>d>>\!\!\S^{2n}H\!\otimes\!\L^{2n,0}E^*
    \!\!\rightarrow\!0 \\[2mm]
    0\!\leftarrow \L^{2n,2n}E^*\!\!@<{\delta}<<\!\! H\!\otimes\!
    \L^{2n,2n-1}E^*\!\!@<{\delta}<<\!\!\! \ldots\!\!\! @<{\delta}<<
    \!\!\S^{2n}H\!\otimes\!\L^{2n,0}E^*\!\!\leftarrow\! 0.
   \end{CD}
  $$
  We note that $\L^{q,0}E^*\cong\L^q E^*$ are canonically isomorphic whereas
  the choice of an isomorphism $\L^{2n,2n-q}E^*\cong \L^qE\cong \L^qE^*$
  amounts to choosing a volume form and a metric respectively on $M$.
  Somewhat more general than the two complexes arising from the de Rham
  complexes are complexes of first order differential operators $d$ and
  $\delta$ first defined by Salamon
  \begin{eqnarray}
    d: \quad L^{-s}\!\otimes\!\S^{k+q}H\!\otimes\!\L^{q,0}E^*\quad
        &\!\!\!\!\longrightarrow\!\!\!\!& L^{-s}\!\otimes\!\S^{k+q+1}H
        \!\otimes\!\L^{q+1,0}E^*\label{dcx}\\[2mm]
    \delta:\;\;L^s\!\otimes\!\S^{k+q}H\!\otimes\!\L^{2n,2n-q}E^*
        &\!\!\!\!\longrightarrow\!\!\!\!& L^s\!\otimes\!\S^{k+q-1}H
        \!\otimes\!\L^{2n,2n-q+1}E^*\label{dscx}
  \end{eqnarray}
  for all even $k\geq 0$ with $s:=\frac{k}{2n+2}$, for odd $k$ the bundles
  involved are in general ill--defined. The twist with the auxiliary line
  bundles $L^s$ and $L^{-s}$ is inserted to make the definition of the
  operators $d$ and $\delta$ independent of the choice of a torsion free
  connection and can be ignored for any other purpose. On the quaternionic
  projective space ${\bf P}^n\BH$ the two complexes above arise from the
  Bernstein--Gelfand--Gelfand resolution of the irreducible representation
  $\S^k(H\oplus E)$ of $\mathbf{PGL}_{n+1}(\BH)$. In this sense the two
  complexes above are curved analogues of the Bernstein--Gelfand--Gelfand
  resolution \cite{barker}.

  Leung and Yi studied the case $k=0$ arising from the de Rham complex and
  proposed to choose a Riemannian metric adapted to the quaternionic structure
  in order to construct an isomorphism $\gamma$ between these two complexes:
  $$
   \begin{CD}
    \BC @>d>> H\otimes\L^{1,0}E^* @>d>>\!\ldots\!@>d>>
    \S^{2n}H\otimes\L^{2n,0}E^* \\
    @V{\gamma}VV @V{\gamma}VV @V{\gamma}VV @V{\gamma}VV \\
    \L^{2n,2n}E^* @<{\delta}<< H\otimes\L^{2n,2n-1}E^*
    @<{\delta}<<\!\ldots\!@<{\delta}<<\S^{2n}H\otimes\L^{2n,0}E^*
    \end{CD}
  $$
  Using this isomorphism they first defined the elliptic second order
  differential operator $\Delta:=(d+\gamma^{-1}\delta\gamma)^2$ and then
  quaternionic analytic torsion as the torsion associated to this Laplacian.
  Note that $\gamma^{-1}\delta\gamma$ will never be the formal adjoint of
  $d$ unless the isomorphism $\gamma$ is parallel. However even if the
  isomorphism $\gamma$ can be chosen to be parallel there remains the
  delicate problem as to its proper choice and the naive choice is
  certainly not the optimal one.

  In order to analyze this problem we will restrict attention to quaternionic
  K\"ahler manifolds or in dimension 4 to half conformally flat Einstein
  manifolds. Recall that choosing a quaternionic K\"ahler metric is equivalent
  to choosing a positive, real section $\sigma_E$ parallel for some torsion
  free connection compatible with the quaternionic structure. Evidently
  its highest power $\frac1{n!}\sigma_E^n$ defines a parallel trivialization
  of all the bundles $L^s,\;s\in\BR$. Any natural choice for $\gamma$ is
  parallel, too, and for appropriate choices of the Hermitian metrics on
  the bundles involved the operator $\gamma^{-1}\delta\gamma$ will be the
  formal adjoint of $d$ as expected. Recall that the operators $d$ and
  $\delta$ for $k=0$ arise as quotient or subcomplexes of the de Rham
  complex on forms. In particular both $d$ and $\delta$ are determined
  by their symbols $\sigma_d[\alpha\otimes\eta]:\;\S^qH\otimes\L^qE^*\,
  \longrightarrow\,\S^{q+1}H\otimes\L^{q+1}E^*$ and $\sigma_\delta[\alpha
  \otimes\eta]:\;\S^qH\otimes\L^{2n,\,2n-q}E^*\,\longrightarrow\,\S^{q-1}H
  \otimes\L^{2n,\,2n-q+1}E^*$ respectively, which are given by
  \begin{equation}\label{keqz}
   \sigma_d[\alpha\otimes\eta]\;:=\;
   \tfrac1{q+1}\alpha^\flat\cdot\otimes\eta\wedge\qquad\quad
   \sigma_\delta[\alpha\otimes\eta]\;:=\;
   (-1)^q\alpha\i\otimes\mathrm{id}\otimes\eta\wedge
  \end{equation}
  according to formula (\ref{wf}). In fact they are the composition of their
  symbol with the covariant derivate with respect to some adapted torsion
  free connection. Similarly the operators $d$ and $\delta$ are defined simply
  by specifying their symbols
  $\sigma_d[\alpha\otimes\eta]:\,L^{-s}\otimes\S^{k+q}H
  \otimes\L^qE^*\longrightarrow L^{-s}\otimes
  \S^{k+q+1}H\otimes\L^{q+1}E^*$ and $\sigma_{\delta}[\alpha\otimes\eta]:
  L^s\!\otimes\S^{k+q}H\otimes\L^{2n,2n-q}E^*\!\longrightarrow\!L^s
  \!\otimes \S^{k+q-1}H\otimes\L^{2n,2n-q+1}E^*$ generalizing (\ref{keqz}):
  \begin{equation}\label{symbol}
   \sigma_d[\alpha\otimes\eta]\,:=\,\tfrac1{k+q+1}\mathrm{id}
   \otimes\alpha^\flat\cdot\otimes\eta\wedge\qquad
   \sigma_\delta[\alpha\otimes\eta]\,:=\,(-1)^q\mathrm{id}
   \otimes\alpha\i\otimes\mathrm{id}\otimes\eta\wedge\,.
  \end{equation}

  Given now a Riemannian metric on $M$ adapted to the quaternionic
  structure or equivalently a real, positive section $\sigma_E$ of
  $\L^2E^*$ we may naively choose $\gamma$ to be the musical isomorphism
  $\flat:\;\L^qE^*\longrightarrow\L^qE\cong\L^{2n,2n-q}E^*$. However the
  associated formal Laplacian $(d+\flat^{-1}\delta\flat)^2$ fails to have
  the right symbol to be properly called a Laplacian even for $k=0$, in
  fact its symbol at an isotropic covector $\alpha\otimes\eta\in(H\otimes E)^*$
  does not act trivially on $\S^qH\otimes\L^qE^*,\,q>0$:
  \begin{eqnarray*}
   \sigma_{(d+\flat^{-1}\delta\flat)^2}[\alpha\otimes\eta]&=&
    \{\,\sigma_d[\alpha\otimes\eta],\,\sigma_{\flat^{-1}\delta\flat}
    [\alpha\otimes\eta]\,\}\\
   &=& -\;\tfrac1{q(q+1)}\,\alpha^\flat\cdot\alpha\i\otimes
        \eta^\flat\wedge\eta\i
  \end{eqnarray*}

  One way to understand this problem is to observe that the vector bundle
  $\L^qE\cong\L^qE^*$ involved is no longer irreducible under the holonomy
  group $\Sp(1)\cdot\Sp(n)\subset\Sp(1)\cdot\GL_n(\BH)$ of the Levi--Civita
  connection of a quaternionic K\"ahler manifold but decomposes into parallel
  subbundles according to
  $$
   \L^q\;E^*\quad=\quad\bigoplus_{r=0\atop r\equiv q\,(2)}^{q\wedge(2n-q)}
   \;\L^r_\circ E^*\,,
  $$
  where the trace free exterior power $\L^r_\circ E^*$ is the kernel
  of the contraction with the dual of the symplectic form. Consider now
  the spinor representation of $\Sp(1)\times\Sp(n)$ (cf.~\cite{barker},
  \cite{wg}):
  $$
   \pi_{\raise.1em\hbox to 0pt{\hskip.1em$\scriptstyle\slash$\hss}S}
   \quad=\quad\bigoplus_{r=0}^n\;\pi_{\raise.1em\hbox to 0pt{\hskip.1em
   $\scriptstyle\slash$\hss}S_r}\quad:=\quad
   \bigoplus_{r=0}^n\;\S^{n-r}\pi_H\otimes\L^r_\circ\pi_E^*\,.
  $$
  As noted in \cite{sal} there is a $\BZ_2$--graded isomorphism
  of vector bundles
  \begin{eqnarray*}
   \bigoplus_{\scriptscriptstyle q=0}^{\scriptscriptstyle 2n}\,
   \S^{k+q}H\otimes\L^qE^*
   &\!\cong\!&\bigoplus_{\scriptscriptstyle q=0}^{\scriptscriptstyle 2n}\!
          \bigoplus_{r=0\atop r\equiv q\,(2)}^{\scriptscriptstyle q
          \wedge(2n-q)}
          \;\S^{k+q}H\otimes\L^r_\circ E^*\\
   &\!\cong\!&\bigoplus_{\scriptscriptstyle r=0}^{\scriptscriptstyle n}\!
          \bigoplus_{q=r\atop q\equiv r\,(2)}^{\scriptscriptstyle 2n-r}
          \;\S^{k+q}H\otimes\L^r_\circ E^*
   \;\;\cong\;\;\SB\otimes\S^{k+n}H
  \end{eqnarray*}
  and it is natural to ask whether the isomorphism $\gamma$ we are looking
  for can be chosen in such a way that the operator $d\,+\,\gamma^{-1}\delta
  \gamma$ and the twisted Dirac operator on $\SB\otimes\S^{k+n}H$ are
  intertwined. A complete answer to that question involves the following
  technical lemma:

  \begin{lemma}\label{cm}
   Consider the subspace $\S^{k+q}H\otimes\L^r_\circ E^*$ of $\S^{k+q}H
   \otimes\L^qE^*$ and set $l:=\frac{q-r}2$ for convenience. The Clifford
   module structure of the twisted spinor bundle $\bigoplus_q\,\S^{k+q}H
   \otimes\L^qE^*\;\cong\;\SB\otimes\S^{k+n}H$ gives rise to the following
   Clifford multiplication on this subspace:
   \begin{eqnarray}\label{cmult}
    \tfrac1{\sqrt{2}}\,(h\otimes e)\bullet&=&\;\;
    \tfrac1{k+q+1}\,h\cdot\otimes e^\sharp\wedge_\circ
    \quad\;-\;\;\;\tfrac1{k+q+1}\,\tfrac{l+1}{n-r+1}\,h\cdot\otimes e\i\\[1mm]
    && +\tfrac{k+n+l+1}{k+q+1}\,h^\sharp\i\otimes e^\sharp\wedge_\circ
    \;+\;\tfrac{k+q-l}{k+q+1}\,\tfrac{n-r-l+1}{n-r+1}\,h^\sharp\i\otimes e\i\,.
    \nonumber
   \end{eqnarray}
  \end{lemma}

  Naturally the Clifford multiplication is defined only up to conjugation
  by an $\Sp(1)\cdot\Sp(n)$--equivariant isomorphism and this freedom allows
  to chose the first two constants more or less arbitrarily as long as a
  simple compatibility condition is met; the other two constants are fixed
  uniquely by this choice. For the time being we state the formula as it is
  with constants convenient to relate the operator $d+\gamma^{-1}\delta\gamma$
  to a twisted Dirac. Of course it is only natural to be curious about a
  satisfactory explanation for the constants appearing in this formula
  (\ref{cmult}).

  Checking the Clifford relation for the Clifford multiplication (\ref{cmult})
  directly seems prohibitively difficult. However it is much easier to see
  that the anticommutator $\{(h\otimes e)\bullet,(\tilde h\otimes\tilde e)
  \bullet\}$ maps $\S^{k+q}H\otimes\L^r_\circ E$ to itself. We will give a
  brief sketch of this calculation before we proceed to the actual proof of
  Lemma \ref{cm} to vindicate formula (\ref{cmult}) and to convince the
  reader that the constants above are much less arbitrary as they may seem
  at first glance.

  Consider the components of the anticommutator $\{(h\otimes e)\bullet,
  (\tilde h\otimes\tilde e)\bullet\}$ mapping $\S^{k+q}H\otimes\L^r_\circ E^*$
  to the various summands of $\SB\otimes\S^{k+n}H$. By definition $e^\sharp
  \wedge_\circ$ is the composition of $e\wedge$ with the projection to the
  trace free subspace $\L^{r+1}_\circ E^*$ of $\L^{r+1}E^*$, in particular
  $e^\sharp\wedge_\circ$ and $\tilde e^\sharp\wedge_\circ$ anticommute as
  do $e^\sharp\wedge$ and $\tilde e^\sharp\wedge$. Hence the components of
  the anticommutator mapping to $\S^{k+q\pm2}H\otimes\L^{r\pm2}_\circ E^*$
  certainly vanish. Moreover there is a fundamental identity on two dimensional
  symplectic vector spaces like $H$, namely $\sigma_H(h,a)\tilde h\,-\,\sigma_H
  (\tilde h,a)h\;=\;\sigma_H(h,\tilde h)a$ for all $h,\tilde h$ and $a\in H$,
  which implies the identity $h\cdot\tilde h^\sharp\i\,-\,\tilde h\cdot
  h^\sharp\i\;=\;(k+q)\sigma_H (h,\tilde h)$ on $\S^{k+q}H$ or:
  $$
   (k+q+2)\,h\cdot\tilde h^\sharp\i\;+\;(k+q)\,h^\sharp\i\tilde h\cdot
   \;\;=\;\;
   (k+q+1)\,\left(\;h\cdot\tilde h^\sharp\i\;+
   \;\tilde h\cdot h^\sharp\i\;\right)\,.
  $$
  Using this identity the component of the composition $(h\otimes e)\bullet
  (\tilde h\otimes\tilde e)\bullet$ mapping $\S^{k+q}H\otimes\L^r_\circ E^*$
  to $\S^{k+q}H\otimes\L^{r+2}_\circ E^*$ can be written
  \begin{eqnarray*}
  \lefteqn{\tfrac{k+n+l+1}{k+q+1}\,\left(\;\tfrac1{k+q}\,h\cdot\tilde h^\sharp\i
   \otimes e^\sharp\wedge_\circ\tilde e^\sharp\wedge_\circ\;+\;
   \tfrac1{k+q+2}\,h^\sharp\i\tilde h\cdot\otimes e^\sharp\wedge_\circ
   \tilde e^\sharp\wedge_\circ\;\right)}\qquad\qquad\qquad\\
   \qquad\qquad\qquad&=&\tfrac{k+n+l+1}{(k+q)(k+q+2)}\,\left(\,h\cdot\tilde
    h^\sharp\i\;+\;\tilde h\cdot h^\sharp\i\;\right)\otimes e^\sharp
    \wedge_\circ\tilde e^\sharp\wedge_\circ\,,
  \end{eqnarray*}
  which is skew in $h\otimes e$ and $\tilde h\otimes\tilde e$ and hence
  vanishes upon symmetrization. The same argument with a different leading
  constant shows that the anticommutator does not map $\S^{k+q}H\otimes
  \L^r_\circ E^*$ to $\S^{k+q}H\otimes\L^{r-2}_\circ E^*$ either. Completely
  analogous arguments replacing the fundamental identity of two dimensional
  symplectic vector spaces by $e\i\tilde e^\sharp\wedge_\circ\,+\,\tilde
  e^\sharp\wedge_\circ e\i\;=\;\sigma_E(e,\tilde e)\,+\,\frac1{n-r+1}e^\sharp
  \wedge_\circ \tilde e\i$ on $\L^r E^*$ (cf.~\cite{ksw}), more usefully
  written as
  \begin{eqnarray*}
   \lefteqn{(n-r+1)\,e\i\tilde e^\sharp\wedge_\circ\;+\;(n-r)\,e^\sharp
    \wedge_\circ\tilde e\i}\qquad\qquad\\
   \qquad\qquad&=&\;(n-r+1)\,\left(\;\sigma_E(e,\tilde e)\;+\; e^\sharp
   \wedge_\circ \tilde e\i\;-\;\tilde e^\sharp\wedge_\circ e\i\;\right)
  \end{eqnarray*}
  show that the components of the anticommutator $\{(h\otimes e)\bullet,
  (\tilde h\otimes\tilde e)\bullet\}$ mapping to $\S^{k+q\pm2}H\otimes
  \L^r_\circ E^*$ vanish, too, consequently the anticommutator maps
  $\S^{k+q}H \otimes\L^r_\circ E^*$ to itself as claimed.

  \beginProof
  Let us choose embeddings $\iota_{q,\,r}:\;\S^{k+q}H\,\longrightarrow
  \,\S^{n-r}H\otimes\S^{k+n}$H for all $r\leq q\leq 2n-r$ with $q\equiv
  r\;(2)$, which piece together to an isomorphism:
  $$
   \iota\;:=\;\oplus(\iota_{q,\,r}\otimes\textrm{id}):
   \quad\bigoplus_{{q\equiv r\,(2)\atop r\leq q\leq 2n-r}}\,\S^{k+q}H
   \otimes\L^r_\circ E^*\;\;\;\longrightarrow\;\;\;\SB\otimes\S^{k+n}H\,.
  $$
  The diagonal multiplication $\sigma\cdot:\,\S^sH\otimes\S^tH\longrightarrow
  \S^{s+1}H\otimes\S^{t+1}H$ with $\sigma$ and the Pl\"ucker map $\textrm{Pl}:
  \S^sH\!\otimes\S^tH\longrightarrow\S^{s+1}H\!\otimes\S^{t-1}H$ give in fact
  rise to an embedding
  $$
   \iota_{q,\,r}\quad:=\quad{1\over l!}\textrm{Pl}^l\;{1\over(n-r-l)!}\;
   (\sigma\cdot)^{n-r-l}
  $$
  with $l:=\frac{q-r}2$ and $\S^{k+q}H\,\cong\,\BC\otimes\S^{k+q}H$. To
  make sense out of this expression we need to choose a pair $\{h_\mu\},
  \,\{h^\vee_\nu\}$ of dual bases for $H$ and $H^*$ to fix $\sigma\cdot:=
  \sum(h_\nu^\vee)^\flat\cdot\otimes h_\nu\cdot$ and $\textrm{Pl}:=\sum
  h_\nu\cdot\otimes h_\nu^\vee\i$ explicitly. It is not difficult to check
  that in terms of the embeddings $\iota_{q,\,r}$ the symmetric product with
  $h\in H$ or the contraction with $\alpha\in H^*$ in the first factor of
  $\S^{n-r}H\otimes\S^{k+n}H$ is expressed by the following formulas:
  \begin{eqnarray*}
   (h\cdot\otimes 1)\,(\iota_{q,\,r}\omega) &=&
   \;\;\;\;\tfrac{l+1}{k+q+1}\,\iota_{q+1,\,r-1}(h\cdot\omega)\;+\;
   \tfrac{n-r-l+1}{k+q+1}\,\iota_{q-1,\,r-1}(h^\sharp\i\omega)\\[.1ex]
   (\alpha\i\otimes 1)\,(\iota_{q,\,r}\omega) &=&
   -\,\tfrac{k+q+1-l}{k+q+1}\,\iota_{q+1,\,r+1}(\alpha^\flat\cdot\omega)\,+\,
   \tfrac{k+n+l+1}{k+q+1}\,\iota_{q-1,\,r+1}(\alpha\i\omega)\,.
  \end{eqnarray*}
  According to the formula for the Clifford multiplication in the
  untwisted case $\SB\cong\bigoplus_r\S^{n-r}H\otimes\L^r_\circ E^*$
  given in \cite{ksw} the twisted Clifford multiplication on
  $\S^{k+q}H\otimes\L^r_\circ E^*\,\subset\,\SB\otimes\S^{k+n}H$
  becomes
  \begin{eqnarray*}
   \tfrac1{\sqrt{2}}(h\otimes e)\bullet
   &=& \iota^{-1}\circ\left(h\cdot\otimes e\i\otimes\mathrm{id}
       \;-\;\tfrac1{n-r}h^\sharp\i\otimes e^\sharp\wedge_\circ\otimes
       \mathrm{id}\right)\circ\iota\\
   &=& \tfrac{l+1}{k+q+1}\,h\cdot\otimes e\i\;\;+\;\;\tfrac{n-r-l+1}{k+q+1}\,
       h^\sharp\i\otimes e\i\\[1mm]
   &&  \quad+\;\tfrac1{n-r}\,\tfrac{k+q+1-l}{k+q+1}\,h\cdot\otimes e^\sharp
       \wedge_\circ\;-\;\tfrac1{n-r}\,\tfrac{k+n+l+1}{k+q+1}\,h^\sharp\i
       \otimes e^\sharp\wedge_\circ
  \end{eqnarray*}
  under the isomorphism $\iota$, the change of sign in the first line is
  due to the fact that we are working with $\L^r_\circ E^*$ instead of
  $\L^r_\circ E$. Evidently this Clifford multiplication is conjugated
  to the Clifford multiplication stated in (\ref{cmult}) under the
  $\Sp(1)\cdot\Sp(n)$--equivariant isomorphism of $\bigoplus\S^{k+q}H
  \otimes\L^r_\circ E^*$, which is $(-1)^l(k+q-l)!(n-r)!$ on
  $\S^{k+q}H\otimes\L^r_\circ E^*$.
  \endProof

  The straightforward embeddings $\L^r_\circ E^*\longrightarrow\L^q E^*$
  and $\L^r_\circ E^*\longrightarrow\L^{2n,2n-q}E^*$ sending $\psi\in\L^r_\circ
  E^*$ to $\frac1{l!}(\sigma_E\wedge)^l\psi$ and $\frac1{n!}\sigma_E^n\otimes
  \frac1{(n-r-l)!}(\sigma_E\wedge)^{n-r-l}\psi$ respectively translate the
  symbols of the operators $d$ and $\delta$ given explicitly in (\ref{symbol})
  into the following maps on the subspace $\S^{k+q}H\otimes\L^r_\circ E^*$ of
  $\S^{k+q}H\otimes\L^qE^*$ and $\S^{k+q}H\otimes\L^{2n,2n-q}E^*$:
  \begin{eqnarray*}
   \sigma_d[h^\sharp\otimes\eta]&=&
   \tfrac1{k+q+1}\,h\cdot\otimes\eta\wedge_\circ\;-\;
   \tfrac1{k+q+1}\,\tfrac{l+1}{n-r+1}\,h\cdot\otimes\eta^\flat\i\\
   (-1)^q\sigma_\delta[h^\sharp\otimes\eta]&=&
   \quad\quad h^\sharp\i\otimes\eta\wedge_\circ\;-\;\tfrac{n-r-l+1}{n-r+1}\,
   h^\sharp\i\otimes\eta^\flat\i
  \end{eqnarray*}
  Comparing this with the formula (\ref{cmult}) for the Clifford
  multiplication we immediately deduce the following proposition
  which is the main result of this section:

  \begin{prop}\label{iso}
   Identify the spaces $\S^{k+q}H\otimes\L^r_\circ E^*,\,r\leq q\leq 2n-r$
   with subspaces of both $\S^{k+q}H\otimes\L^qE^*$ and $\S^{k+q}H\otimes
   \L^{2n,2n-q}E^*$ as above and consider the $\Sp(1)\cdot\Sp(n)$--equivariant
   isomorphism
   $$
    \gamma:\quad\;\S^{k+q}H\otimes\L^qE^*\;\longrightarrow
    \;\S^{k+q}H\otimes\L^{2n,2n-q}E^*\,,
   $$
   which is $(-1)^l\frac{(k+q-l)!(k+n+l+1)!}{(k+q+1)!}$ on these subspaces.
   The differential operator
   $$
    D_{\S^{k+n}H}\quad:=\quad\sqrt{2}\,(d+\gamma^{-1}\circ\delta\circ\gamma)
   $$
   is the twisted Dirac operator on $\SB\otimes\S^{k+n}H\,\cong\,\bigoplus_q
   \S^{k+q}H\otimes\L^qE^*$.
  \end{prop}

  Moreover $\gamma$ is uniquely characterized by this property up to an
  overall constant with respect to the present choice of the Clifford
  multiplication (\ref{cmult}) and the symbols (\ref{symbol}) of the
  operators $d$ and $\delta$. Other conventions simply conjugate $d,
  \,\delta$ and $D_{\S^{k+n}H}$ by $\Sp(1)\cdot\GL_n(\BH)$-- and $\Sp(1)\cdot
  \Sp(n)$--equivariant isomorphisms respectively leading essentially to
  the same conclusion but with appropriately conjugated $\gamma$. It
  is important however to note that the operator $D^2_{\S^{k+n}H}$
  respects the decomposition of $\SB\otimes\S^{k+n}H$ into $\Sp(1)
  \cdot\Sp(n)$--irreducible subspaces and is hence genuinely defined
  independent of all choices.

  Since the operator $d+\gamma^{-1}\delta\gamma$ is a twisted Dirac
  operator on a quaternionic K\"ahler manifold the cohomology of the
  complexes (\ref{dcx}) and (\ref{dscx}) can be presented by harmonic
  twisted spinors. Quite a lot is known about the existence of harmonic
  spinors in this situation and consequently about the cohomology of these
  complexes (\cite{naga}, \cite{sw}). In particular the $d$--complex is
  acyclic for all even $k\geq 0$ except in degree $q=0$, if the scalar
  curvature $\kappa>0$ is positive. Moreover it is assumed that its
  cohomology in degree zero governs the classification of quaternionic
  K\"ahler manifolds with $\kappa>0$. If the scalar curvature $\kappa<0$
  is negative, then the $d$--complex is acyclic except in degree $q=2n$
  (sic!) for all even $k>0$, but for $k=0$ it has trivial cohomology
  $\BC$ in degree $q=0$ and it may have exceptional cohomology in degrees
  $q=n,\ldots,2n$. In the hyperk\"ahler case $\kappa=0$ the cohomology of
  the $d$--complex can be represented by holomorphic forms and thus
  faithfully reflects the decomposition of the manifold into irreducible
  factors.

  For our calculations we are also interested in twisted versions of the
  complexes introduced above. However extra curvature terms arising from
  a twisting bundle $\F$ will spoil $d^2=0$ unless the curvature of $\F$
  will be an antiselfdual two form, i.~e.~a section of $\S^2E\otimes
  \mathrm{End}\;\F\subset\L^2(TM\otimes_\BR\BC)\otimes\mathrm{End}\;\F$.
  Consequently we restrict ourselves to Hermitian vector bundles $\F$ with
  an antiselfdual Hermitian connection. Associated to such an antiselfdual
  bundle $\F$ and all even $k\geq 0$ are twisted versions
  $$
   d_\F:\;\S^{k+q}H\otimes\L^{q,\,0}E^*\otimes\F\;\longrightarrow\;
    \S^{k+q+1}H\otimes\L^{q+1,\,0}E^*\otimes\F
  $$
  and
  $$
   \delta_\F:\;\S^{k+q}H\otimes\L^{2n,\,2n-q}E^*\otimes\F\longrightarrow
    \S^{k+q-1}H\otimes\L^{2n,\,2n-q+1}E^*\otimes\F
  $$
  of the elliptic complexes considered above. The cohomology of the
  $d$--complex defines the quaternionic cohomology $H^{*,\,k}(M,\F)$ of $\F$.
  Let $\lpentagon_{q,k}$ denote the operator
$\lpentagon_k:=(d_\F\,+\,\gamma^{-1}
  \delta_\F\gamma)^2$ restricted to $\S^{k+q}H\otimes\L^{q,\,0}E^*\otimes\F$
  with spectrum $\sigma(\lpentagon_{q,k})$. The usual arguments of Hodge theory
  imply that quaternionic cohomology can be represented by harmonic sections
  $H^{q,k}(M,\F)\cong\mathrm{ker}\,\lpentagon_{q,k}$. In the quaternionic
  K\"ahler case the twisted complexes are related to the Dolbeault complex
  of suitable holomorphic vector bundles on the twistor space via the
  Penrose transform, in particular the main motivation for studying
  these complexes arise from complex geometry (\cite{naga},\cite{mamo}).

  Consider an isometry $g$ of the quaternionic K\"ahler manifold $M$,
  preserving the quaternionic structure (e.~g.~the identity). Assume
  furthermore an isometry of vector bundles $g^\F:\F\to g^*\F$.
  Then the quaternionic torsion is defined via the zeta function
  $$
   Z_g(s):=\sum_{q=0}^{2n} (-1)^{q+1}q \sum_{\lambda\in\sigma(\spentagon_{q,k})
   \atop\lambda\neq0}\lambda^{-s}
   \Tr g^*_{|{\rm Eig}_\lambda(\spentagon_{q,k})}
  $$
  for Re $s\gg 0$. This zeta function has a meromorphic continuation to the
  complex plane which is holomorphic at $s=0$ by a general result by Donnelly
  (\cite{Donn}).

  \begin{defin}
   The equivariant quaternionic analytic torsion is defined as
   $$
    T^k_g(M,\F):=Z'_g(0)\,\,.
   $$
  \end{defin}

  Similarly, one can define an equivariant Quillen metric on the equivariant
  determinant of the quaternionic cohomology. Let $g$ be an isometry 
of an hermitian vector space
$E$. Let $\Theta$ denote  the set of eigenvalues $\zeta$ of $g$ with
associated eigenspaces $E_\zeta$. The
$g$-{\bf equivariant determinant} of $E$ is defined as
$$ {\det}_g E:=\bigoplus_{\zeta\in\Theta} \det E_\zeta.
$$ The $g$-{\bf equivariant metric} associated to the metric on $E$ is the map
\begin{eqnarray*}
\log\|\cdot\|^2_{{\det}_g E}:{\det}_g E&\to&{\bf C}\\
(s_\zeta)_\zeta&\mapsto&\sum_{\zeta\in\Theta}
\log\|s_\zeta\|^2_\zeta\cdot\zeta,
\end{eqnarray*} where $\|\cdot\|^2_\zeta$ denotes the induced metric on
$\det E_\zeta$. Now in our situation the isometry $g$ induces an isometry
$g^*$ of the cohomology
$H^{q,k}(M,\F)\cong\ker
\lpentagon_{q,k}$ equipped with the restriction of the $L^2$-metric.
\begin{defin} Set $\lambda_g(M,\F):=\big[{\det}_g H^{q,k}(M,\F)\big]^{-1}$.
  The equivariant Quillen metric on
$\lambda_g(M,E)$ is defined as
\begin{equation}
\log\|\cdot\|^2_{Q,\lambda_g(M,E)}:=\log\|\cdot\|^2_{L^2,\lambda_g(M,\F)}
-T^k_g(M,\F).
\end{equation}
\end{defin}

\section{Quaternionic torsion for symmetric spaces}

  On a symmetric space $G/K$ Partharasarty's formula relates the squares of
  twisted Dirac operators to the Casimirs of $G$ and $K$. In consequence the
  operators $2\,(d_\F\,+\,\gamma^{-1}\delta_\F\gamma)^2$, which are squares of
  twisted Dirac operators on every quaternionic K\"ahler manifold, can be
  expressed in terms of the Casimirs of $G$ and $K$ on every quaternionic
  K\"ahler symmetric space. Recall that these two Casimirs induce an elliptic
  second order differential operator and a curvature operator respectively on
  every homogeneous vector bundle on $G/K$. Analogues of these two operators
  are defined in \cite{sw} for all vector bundles associated to the holonomy
  bundle $\mathrm{Hol}(M)$ of an arbitrary Riemannian manifold $M$ via
  representations of the holonomy group.

  In fact the Levi--Civita connection $\nabla$ of $M$ defines an elliptic
  second order differential operator on every homogeneous vector bundle,
  namely the horizontal Laplacian $\nabla^*\nabla$. On a symmetric space
  $G/K$ with metric induced by the Killing form $B$ of $G$ on $\g=\k\oplus
  \p$ the horizontal Laplacian is the ``Casimir operator'' of $\p$ up to
  sign. It is more difficult to write down the analogue of the Casimir
  of $K$. Consider for this purpose a point $p$ in a Riemannian manifold
  $M$ with holonomy group $\mathrm{Hol}_pM\subset\O(T_pM)$ and holonomy
  algebra $\mathfrak{hol}_pM\subset\L^2(T_pM)$. The completely contravariant
  curvature tensor $R$ of $M$ at $p$ is by its very  definition an element
  of the space $\S^2\mathfrak{hol}_pM\subset\S^2\L^2T_p^*M$ and thus the
  quantization map
  $$
   q:\quad\S\;\mathfrak{hol}_pM\;\longrightarrow\;\mathcal{U}\;
   \mathfrak{hol}_pM,\quad \mathfrak{X}^l\;\longmapsto\;\mathfrak{X}^l
  $$
  defines a curvature term $2\,q(R)\,\in\,\mathcal{U}\,\mathfrak{hol}_pM$
  acting on every vector bundle on $M$ associated to the holonomy bundle
  $\mathrm{Hol}(M)$. Straightforward computation shows that this curvature
  term reduces to the Casimir operator of $K$ with respect to the restriction
  of the Killing form $B$ to $\k$ on every homogeneous vector bundle on $G/K$.
  Consequently the elliptic differential operator
  $$
   \Delta \quad:=\quad \nabla^*\nabla\;\;+\;\;2\,q(R)
  $$
  agrees with the Casimir operator of $G$ on the symmetric space $M\;=\;G/K$.
  The operator $\Delta$ allows us to write the Bochner--Weitzenb\"ock formula
  for a twisted Dirac operator $D_{\mathcal{R}}$ on a twisted spinor bundle
  $\SB\otimes\mathcal{R}$ associated to the holonomy bundle in the form:

  \begin{equation}\label{twis}
   D^2_{\mathcal{R}}\quad=\quad\Delta\;\;+\;\;\frac\kappa8\;\;-\;\;
   \mathrm{id}_\SB\,\otimes\,2\,q(R)
  \end{equation}

  \noindent
  where $\kappa$ is the scalar curvature of $M$ (cf.~\cite{sw}). We will
  employ this formula for tensor products $\mathcal{R}=\S^{k+n}H\otimes\F$
  of $\S^{k+n}H$ with antiselfdual homogeneous vector bundles $\F$ on a
  quaternionic symmetric space $G/K$. Its isotropy group $K=\Sp(1)\cdot
  K':=(\Sp(1)\times K')/\BZ_2$ splits almost into a direct product and by
  definition an antiselfdual homogeneous vector bundle $\F$ is associated
  to a representation on which $\Sp(1)\subset K$ acts trivially. Hence $\F$
  is induced by a representation of $K_\circ:=K'/\BZ_2$. In this case
  formula (\ref{twis}) provides the following corollary to Proposition
  \ref{iso}:

  \begin{cor}\label{Cas}
   Let $\F$ be an antiselfdual homogeneous vector bundle on a quaternionic
   K\"ahler symmetric space $G/K$, i.~e.~the subgroup $\Sp(1)\subset K$
   acts trivially on the corresponding representation of $K$. The square
   of the operator $d_\F\,+\,\gamma^{-1}\delta_\F\gamma$ for even $k\geq 0$
   can be expressed as:
   $$
    (\,d_\F\,+\,\gamma^{-1}\,\delta_\F\,\gamma\,)^2\quad=\quad
    {1\over2}\left(\,\mathrm{Cas}_G\,-\,\frac\kappa8\,\frac{k(k+2n+2)}{n(n+2)}
    \,-\,\mathrm{Cas}^\F_K\,\right)\,.
   $$
  \end{cor}

  \beginProof
  The algebraic relation between $d\,+\,\gamma^{-1}\,\delta\,\gamma$
  and $D_{\S^{k+n}H}$ proved in Proposition \ref{iso} remains valid
  under arbitrary twists. In particular we may use the identification
  of $\SB\otimes\S^{k+n}H\otimes\F$ with $\bigoplus_q\S^{k+q}H
  \otimes\L^qE^*\otimes\F$ to write the operator $(d_\F\,+\,\gamma^{-1}\,
  \delta_\F\gamma)$ as a twisted Dirac operator:
  $$
   \sqrt{2}\,(\,d_\F\,+\,\gamma^{-1}\,\delta_\F\,\gamma\,)\quad=\quad
   D_{\S^{k+n}\otimes\F}
  $$
  Equation (\ref{twis}) relates the operator $D^2_{\S^{k+n}H\otimes\F}$
  to the Casimirs of $G$ and $K$:
  $$
   2\,(\,d_\F\,+\,\gamma^{-1}\,\delta_\F\,\gamma\,)^2
   \quad=\quad
   D^2_{\S^{k+n}H\otimes\F}
   \quad=\quad
   \mathrm{Cas}_G\,+\,{\kappa\over 8}\,-\,
   \mathrm{Cas}_K^{\S^{k+n}H\otimes\F}
  $$
  However the Lie algebra of $K$ splits into commuting subalgebras
  $\k=\mathfrak{sp}(1)\oplus\k_\circ$ and $\mathfrak{sp}(1)$ acts
  trivially on the representation corresponding to $\F$ by assumption
  whereas $\k_\circ$ acts trivially on the representation corresponding
  to $\S^{k+n}H$. Hence the Casimir of $K$ on $\S^{k+n}H\otimes\F$ is
  the sum:
  $$
   \mathrm{Cas}_K^{\S^{k+n}H\otimes\F}
   \quad=\quad\frac\kappa8\,\frac{(k+n)(k+n+2)}{n(n+2)}\;+\;1\otimes
   \mathrm{Cas}_K^\F\,.
  $$
  In fact the Casimir of $\S^{k+n}H$ is proportional to $(k+n)(k+n+2)$
  and necessarily equals $\frac\kappa8$ for $k=0$, because $\SB\otimes\S^nH$
  occurs in the forms.
  \endProof

Fix a maximal torus $T$ of $K=\Sp(1)\cdot K'$ containing a maximal torus of
$\Sp(1)$.
  T is automatically a maximal torus of $G$. As $T$ contains a maximal torus of
  $\Sp(1)$ the subalgebra $\mathfrak{sp}(1)$ is invariant under $T$ and we
denote
  its weights by $-2\alpha,\,0,\,2\alpha$. Under the action of
$\mathfrak {sp}(1)$  the Lie algebra of $\g$ splits into
$\g^{-2\alpha}\oplus\g^{-\alpha}\oplus\g^0\oplus\g^{\alpha}\oplus\g^{2\alpha}$.
We
choose an ordering of the roots of $G$ such that $2\alpha$ is the highest
root and
the weights of $\g^\alpha$ and $\g^{2\alpha}$ are positive.
In particular $\a$ is the positive weight of $H$. We denote the set of
positive roots
by $\Sigma^+$.

Set ${\frak t}_{\rm reg}:=\{X\in{\frak t}|\beta(X)\notin {\bf Z} \
\forall\beta\in\Sigma\}$. For $X\in\frak t$ let $e^X\in T$ denote the
associated group
element. Let $\rho$ denote half the sum of the positive weights of $G$ and
define
similarly $\rho_K$ etc. Let $W_G, W_K$ etc. denote the Weyl groups. Set for
$b\in{\frak t}^*$
$$
{\rm Alt}_G\{b\}:=\sum_{w\in W_G}{\rm sign}(w)e^{2\pi i w b}\,\,.
$$
We denote the
$G$-representation with highest weight $\lambda$ by $V^G_{\rho+\lambda}$
and its
character is denoted by $\chi_{\rho+\lambda}$. In general, for a weight
$\lambda$
and $X\in{\frak t}_{\rm reg}$ we define $\chi_{\rho+\lambda}$ by the Weyl
character
formula
$$
\chi_{\rho+\lambda}(e^X):=\frac{{\rm
Alt}_G\{\rho+\lambda\}(X)}{{\rm Alt}_G\{\rho\}(X)}
$$
with ${\rm Alt}_G\{\rho\}(X)=\prod_{\beta\in\Sigma^+}2i\sin\pi\beta(X)$.
For an irreducible representation $\pi$, we shall denote the sum of $\rho$
and the
highest weight by $b_\pi$. Thus the Casimir acting on $V_\pi$ is given by
$\|b_\pi\|^2-\|\rho\|^2$.

An irreducible $K$-representation $V^K_{\rho_K+\lambda}$
induces a
$G$-invariant vector bundle $\F^K_{\rho_K+\lambda}$ on $M$. As
$V^K_{\rho_K+\lambda}$ carries a $K$-invariant Hermitian metric which is
unique up to a factor, we get corresponding $G$-invariant metrics on
$\F^K_{\rho_K+\lambda}$.
Consider a $K_\circ$-representation
$V^{K_\circ}_{\rho_{K_\circ}+\lambda_\circ}$ of
highest weight
$\lambda_\circ$ and the induced equivariant bundle $\F$ on the quaternionic
K\"ahler symmetric space
$G/K$. Set $\lambda:=\lambda_\circ+k\a$. By Corollary
\ref{Cas}, the zeta function defining the torsion $T^k(M,\F)$ of $\F$ equals
\begin{eqnarray*}
Z(s)&=&\sum_{q=1}^{2n}(-1)^{q+1} q\sum_{\pi\ {\rm irr.}}
\left( \frac{2}{\|b_\pi\|^2-\|\rho+\lambda\|^2} \right)^s
\\&&\cdot
\chi_{b_\pi}
\dim \Hom_K(V_\pi,\Lambda^q E\otimes \S^{k+q} H\otimes
V^{K_\circ}_{\rho_{K_\circ}+\lambda_\circ})\,\,.
\end{eqnarray*}
Let $\Theta^E$
denote the representation of $K$ on $E$ and let $\Psi_0$ denote its weights.
Analogously to \cite[Lemma4]{KC} and \cite[Lemma7]{KR} we show
\begin{lemma}\label{occur}
Let $G/K$ be a $n$-dimensional quaternionic K\"ahler symmetric space. For any
irreducible
$G$-representation $(V_\pi,\pi)$ the sum
$$
\sum_{q=1}^{2n}(-1)^q q \dim \Hom_K(V_\pi,\Lambda^q E\otimes \S^{k+q} H \otimes
V^{K_\circ}_{\rho_{K_\circ}+\lambda_\circ})
$$
equals the sum of $-\chi_{\rho+\lambda+\ell (\a+\b)}$ over those $\ell\in\BN$,
$\b\in\Psi_0$ such that $b_\pi$ is in the $W_G$-orbit of $\rho+\lambda+\ell
(\a+\b)$.
\end{lemma}
\beginProof
Let $\chi^K$ denote the virtual $K$-character
$$
\chi^K:=\sum_{q=1}^{2n} (-1)^q q \chi(\Lambda^q E\otimes \S^{k+q}H)\,\,.
$$
Notice that
$$
\chi(\S^{k+q}H)=\frac{e^{2\pi i(k+q+1)\a}-e^{-2\pi i(k+q+1)\a}}{e^{2\pi
i\a}-e^{-2\pi
i\a}}
$$
and, for $s\in\BR$,
$$
\sum_{q=0}^{2n} (-s)^q \chi(\Lambda^q E)=\det(1-s \th)\,\,.
$$
Hence
\begin{eqnarray*}
\chi_s&:=&\sum_{q=0}^{2n} (-s)^q \chi(\Lambda^q E \otimes \S^{k+q}H)
\\&&
=\frac1{e^{2\pi i\a}-e^{-2\pi i\a}}\left[e^{2\pi i(k+1)\a} \det(1-s e^{2\pi
i\a}\th)-
e^{-2\pi i(k+1)\a} \det(1-s e^{-2\pi i\a}\th)\right]
\end{eqnarray*}
and
\begin{eqnarray*}
\chi^K&=&\frac\partial{\partial s}_{|s=1}\chi_s
=\frac{e^{2\pi i(k+1)\a}}{e^{2\pi i\a}-e^{-2\pi i\a}} \det(1-e^{2\pi i\a}\th)
\Tr\left[(1-e^{-2\pi i\a}(\th)^{-1})^{-1}\right]
\\&&- \frac{e^{-2\pi i(k+1)\a}}{e^{2\pi i\a}-e^{-2\pi i\a}}
\det(1-e^{-2\pi i\a}\th) \Tr\left[(1-e^{2\pi i\a}(\th)^{-1})^{-1}\right]
\\&=&\frac1{e^{2\pi i\a}-e^{-2\pi i\a}}\prod_{\b\in\Psi_0,\beta>0}
\left(e^{\pi i(\a+\b)}-e^{-\pi i(\a+\b)}\right)
\left(e^{\pi i(\a-\b)}-e^{-\pi i(\a-\b)}\right)
\\&&\cdot
\Big[e^{2\pi i(k+n+1)\a} \sum_{\b\in\Psi_0,\beta>0}\big(\frac1{1-e^{2\pi
i(-\a+\b)}}+
\frac1{1-e^{2\pi i(-\a-\b)}}\big)
\\&&
-e^{-2\pi i(k+n+1)\a} \sum_{\b\in\Psi_0,\beta>0}\big(\frac1{1-e^{2\pi
i(\a+\b)}}+
\frac1{1-e^{2\pi i(\a-\b)}}\big)\Big]
\\&=&
  \frac{\overline{\Alt_G\{\rho\}}}{\overline{\Alt_K\{\rho_K\}}}\cdot
  \frac1{e^{2\pi i\a}-e^{-2\pi i\a}}
\\&&\cdot
\Big[e^{2\pi i(k+n+1)\a} \sum_{\b\in\Psi_0,\beta>0}\big(\frac1{1-e^{2\pi
i(-\a+\b)}}+
\frac1{1-e^{2\pi i(-\a-\b)}}\big)
\\&&
-e^{-2\pi i(k+n+1)\a} \sum_{\b\in\Psi_0,\beta>0}\big(\frac1{1-e^{2\pi
i(\a+\b)}}+
\frac1{1-e^{2\pi i(\a-\b)}}\big)\Big]\,\,.
\end{eqnarray*}
Notice that $\Sigma_G^+\setminus\Sigma_K^+=\{\a\pm\b|\b\in\Psi_0,\beta>0\}$ and
$\Sigma_K^+\setminus\Sigma_{K_\circ}^+=\{2\a\}$; also, $W_K=W_{K_\circ}\times
W_{{\Sp}(1)}$. In particular, $\rho-\rho_{K_\circ}=(n+1)\a$. Thus
\begin{eqnarray*}\lefteqn{
\overline{\Alt_K\{\rho_K\}}{\Alt_K\{\rho_K\}}\overline{\chi^K\cdot
\chi^{K_\circ}_{\rho_{K_\circ}+\lambda_\circ}}
}\\&=&
  \Alt_G\{\rho\}\overline{{\Alt_{K_\circ}
   \{\rho_{K_\circ}+\lambda_\circ\}}}\sum_{\b\in\Psi_0\atop w\in
   W_{{\Sp}(1)}} {\rm sign}(w)\frac{e^{-2\pi i w(\rho-\rho_{K_\circ}+k\a)}}
   {1-e^{2\pi i w(\a+\b)}}\,\,.
  \end{eqnarray*}
  As in \cite[eq. (29)]{KC} we obtain for large $N\in\BN$
  \begin{eqnarray*}\lefteqn{
\dim \Hom_K(V_\pi,\Lambda^q E\otimes \S^{k+q} H\otimes
V^{K_\circ}_{\rho_{K_\circ}+\lambda_\circ})}
\\&=&
\frac1{\#W_K}\int_T
\overline{\Alt_K\{\rho_K\}}{\Alt_K\{\rho_K\}}\overline{\chi_K
\cdot
\chi^{K_\circ}_{\rho_{K_\circ}+\lambda_\circ}}
\chi_\pi \,{\rm dvol}_T
\\&=&
\frac{-1}{\#W_K}\int_T
\Alt_G\{\rho\} \overline{{\Alt_{K_\circ}\{\rho_{K_\circ}+\lambda_\circ\}}}
\\&&\cdot
\sum_{\b\in\Psi_0\atop w\in W_{{\Sp}(1)}} {\rm sign}(w)
\frac{e^{-2\pi i w(\rho-\rho_{K_\circ}+k\a)}(e^{-2\pi i
w(\a+\b)}-e^{-2\pi i N w(\a+\b)})}{1-e^{-2\pi i w(\a+\b)}}
\chi_\pi \,{\rm dvol}_T
\\&=&
\frac{-1}{\#W_K}\sum_{\ell=1}^{N-1} \int_T
\Alt_G\{b_\pi\} \overline{{\Alt_{K_\circ}\{\rho_{K_\circ}+\lambda_\circ\}}}
\\&&\cdot
\sum_{\b\in\Psi_0\atop w\in
W_{{\Sp}(1)}} {\rm sign}(w)
e^{-2\pi i w(\rho-\rho_{K_\circ}+k\a+\ell(\a+\b))}
\,{\rm dvol}_T
\\&=&
\frac{-1}{\#W_K}\sum_{\ell=1}^{N-1} \sum_{\b\in\Psi_0}\int_T
\Alt_G\{b_\pi\} \overline{{\Alt_{K}\{\rho+\lambda+\ell(\a+\b)\}}}
  \,{\rm dvol}_T\,\,.
\end{eqnarray*}
This proves the Lemma the same way as in \cite[p. 100]{KC}.
\endProof

  Set $\Psi_0^+:=\{\beta\in\Psi_0|\langle(\a+\beta)^\vee,\rho+\lambda\rangle
  \geq0\}$ and $\Psi_0^-:=\{\beta\in\Psi_0|\langle(\a+\beta)^\vee,\rho+\lambda
  \rangle<0\}$ with $\beta^\vee=2\beta/\|\beta\|^2$. By Lemma \ref{occur} and
  Theorem \ref{Cas} we find that $Z(s)$ is given by the following formula:
  \begin{theor}
  For $G/K$ quaternionic K\"ahler, the zeta function $Z$ equals
  \begin{eqnarray*}
   Z(s)&=&-2^s\sum_{\b\in\Psi^+_0}\sum_{\ell>\langle(\a+\beta)^\vee,
   \rho+\lambda\rangle}
   \frac{\chi_{\rho+\lambda+\ell(\a+\b)}}{
   \langle 2\rho+2\lambda+\ell(\a+\b),k(\a+\b)\rangle^{s}}
   \\
   &&+2^s\sum_{\b\in\Psi^-_0}\sum_{\ell>-\langle(\a+\beta)^\vee,
   \rho+\lambda\rangle}
   \frac{\chi_{\rho+\lambda+\ell(\a+\b)}}{
   \langle
   2\rho+2\lambda+\ell(\a+\b),\ell(\a+\b)\rangle^{s}}
   \,\,.
  \end{eqnarray*}
  \end{theor}
  This is a zeta function of the form considered in \cite[Lemma 8]{KC}.
  It is actually the very same formula as in \cite[Prop. 5.1]{KK} (see
  also \cite[Theorem 5]{KC}), the only difference being that we consider
  a different kind of symmetric space there. Define for $\phi\in{\bf R}$
  and ${\rm Re}\,s>1$
  \begin{equation}
   \zeta_L(s,\phi)=\sum_{\ell=1}^\infty \frac{e^{i\ell\phi}}{k^s}\,\,.
  \end{equation}
  The function $\zeta_L$ has a meromorphic continuation to the complex plane
  in $s$ which is holomorphic for $s\neq 1$. Set $\zeta'_L(s,\phi):=\partial/
  \partial s (\zeta_L(s,\phi))$. Let $P:{\bf Z}\to{\bf C}$ be a function of
  the form
  \begin{equation}
   P(\ell)=\sum_{j=0}^m c_j\ell^{n_j}
   e^{i\ell\phi_j}\label{sternnn}
  \end{equation}
  with $m\in{\bf N}_0$, $n_j\in{\bf N}_0$, $c_j\in{\bf C}$, $\phi_j\in{\bf R}$
  for all $j$. We define $P^\odd(\ell):=(P(\ell)-P(-\ell))/2$. Also we define
  as in \cite[Section 6]{KC}
  \begin{eqnarray}
   \bz P&:=&\sum_{j=0}^m c_j\zeta_L(-n_j,\phi_j),\\
   \bzs P&:=&\sum_{j=0}^m c_j\zeta_L'(-n_j,\phi_j),\\
   \mbox{and}\qquad P^*(p)&:=&-\sum_{j=0\atop\phi_j\equiv 0{\rm\ mod}2\pi}^m
    c_j \frac{p^{n_j+1}}{4(n_j+1)} \sum_{\ell=1}^{n_j}\frac 1 \ell
  \end{eqnarray}
  for $p\in{\bf R}$.

  Then by \cite[Lemma 8]{KC} we get the same formula as in
  \cite[Theorem 5.2]{KK} (compare also \cite[Theorem 9]{KC}):
  \begin{theor}
   Let $G/K$ be a quaternionic K\"ahler symmetric space. The equivariant
   analytic torsion of $\overline\F$ on $G/K$ is given by
   \begin{eqnarray*}
    \lefteqn{
     T^k(G/K,\overline \F)=-2\sum_{\beta\in\Psi_0} \bzs
     \chi^\odd_{\rho+\lambda-\ell(\a+\beta)}
     -2\sum_{\beta\in\Psi_0}\chi^*_{\rho+\lambda-\ell(\a+\beta)}
     (\langle(\a+\beta)^\vee,\rho+\lambda\rangle)}\\
    &&-\sum_{\beta\in\Psi_0} \bz
    \chi_{\rho+\lambda-\ell(\a+\beta)}\cdot\log\frac2{\|\a\|^2+\|\beta\|^2}-
    \chi_{\rho+\lambda}\sum_{\beta\in\Psi_0^+}\log\frac2{\|\a\|^2+\|\beta\|^2}\\
    &&-\sum_{\beta\in\Psi_0^+}\sum_{\ell=1}^{\langle(\a+\beta)^\vee,
    \rho+\lambda\rangle}\chi_{\rho+\lambda-\ell(\a+\beta)}\cdot\log \ell
    +\sum_{\beta\in\Psi_0^-}\sum_{\ell=1}^{\langle-(\a+\beta)^\vee,
    \rho+\lambda\rangle}\chi_{\rho+\lambda+\ell(\a+\beta)}\cdot\log \ell\,\,.
   \end{eqnarray*}
  \end{theor}

\section{The hyperk\"ahler case}

  Assume that $M$ is a hyperk\"ahler manifold. Then choosing a
  subordinate complex structure is equivalent to fixing isomorphisms
  $T^{(1,0)}M=T^{(0,1)}M=E$ and $H={\cal O}\oplus{\cal O}$. The
  antiselfdual two forms on a hyperk\"ahler manifold can be identified
  with the forms of bidegree $(1,1)$ for any subordinated complex structure,
  so that antiselfdual vector bundles are exactly the vector bundles which
  are holomorphic with respect all these complex structures. In this case
  the twisted $d$--complex reads
  $$
   \cdots\stackrel{\bar\partial}\to
   \Lambda^qT^{*(0,1)}M\otimes{\cal O}^{\oplus(k+q+1)}\otimes\F
   \stackrel{\bar\partial}\to\cdots\,\,,
  $$
  i.~e.~it is essentially equivalent to the twisted Dolbeault complex up to
  the trivial factor $\S^{k+q} H\cong{\cal O}^{\oplus(k+q+1)}$, which is
  needed to make this complex independent of the choice of a complex structure.

Let $\square_q:=(\bar\partial+\bar\partial^*)^2$ denote the Kodaira-Laplace
operator
acting on $$\Gamma^\infty(M,\Lambda^qT^{*(0,1)}M \otimes \F)\,\,.$$ Let
$P^\perp$ be the
projection of this space to the orthogonal complement of ker $\square_q$
and define
$$
\zeta_q(s):=\Tr (\square_q^{-s}P^\perp)
$$
for Re $s\gg0$. Then the quaternionic torsion equals
\begin{equation}
T^k(M,\F)=\sum_{q=0}^{2n} (-1)^{q+1} q(q+k+1) \zeta'_q(0)
=T^0(M,\F)+kT_{\bar\partial}(M,\F)
\end{equation}
with $T_{\bar\partial}$ denoting the holomorphic torsion.
Now let $K_X:=\Lambda^m T^*X$ denote
the canonical line bundle on an $m$-dimensional compact K\"ahler manifold
$X$. Let
  $\F$ denote a holomorphic Hermitian
vector bundle on
$X$. In
\cite[Th. 1.4]{GSZ} it was shown that $\zeta'_q(0)=\zeta'_{m-q}(0)$ and thus
$$
T_{\bar\partial}(X,\F)=(-1)^{m+1}T_{\bar\partial}(X,\F^*\otimes
K_X)\,\,.
$$
In particular,  if $X$ is even-dimensional and spin and $K_X^{1/2}$ denotes a
chosen square root of $K_X$, then $T_{\bar\partial}(X,K_X^{1/2})$ vanishes.
This statement takes a particularly nice form if $ K_X\cong{\cal O}$ as
holomorphic Hermitian bundles, i.e. for Calabi-Yau manifolds equipped with
the K\"ahler-Einstein metric. This has been noticed first by J.-B. Bost and
J.-M. Bismut.

In our case, this implies that for any bundle $\F$ with $\F=\F^*$, $T^k(M,\F)$
is independent of $k$. In particular this holds for $\F={\cal O}$. If furthermore
$n=1$ we find that the torsion equals $T_k(M,\F)=2\zeta'_0(0)$, twice the
determinant on the space of functions. For any holomorphic Hermitian
bundle
$\F$
\begin{eqnarray*}
T^k(M,\F)&=&\sum_{q=0}^{2n} (-1)^{2n-q+1} (2n-q)(2n-q+k+1)
\zeta'_q(0)(M,\F^*)\\ &=&\sum_{q=0}^{2n} (-1)^{q+1} q(q-(4n+k+1))
\zeta'_q(0)(M,\F^*)
\end{eqnarray*}
which can be interpreted as "$T^k(M,\F)=T^{-4n-k-1}(M,\F^*)$" (we did not
define the latter). The quadratic term in $q$ thus does not provide vanishing
results for the quaternionic torsion in contrast to the holomorphic case.

{\bf Example:} Consider a $2n$-dimensional Hermitian vector
space $V$ and a lattice
$\Lambda\subset V$ of maximal rank. Let $\Lambda^\vee$ denote the dual
of $\Lambda$ and let $M:=V/\Lambda$ be the associated flat torus. Then, as
holomorphic
Hermitian bundles,
$$
\Lambda^qT^{*(0,1)}M\cong{\cal O}^{\oplus{2n\choose q}}
$$
and the zeta function defining $T^k(M,{\cal O})$ equals
\begin{eqnarray*}
\sum_{q=0}^{2n} (-1)^qq(q+k+1) \Tr (\square_q^{-s}P^\perp)
&=&\frac{\partial^2}{\partial x^2}_{|x=1}\left(
\sum_{q=0}^{2n} (-1)^q x^{q+1} {2n\choose q}
\right)\cdot
\Tr (\square_q^{-s}P^\perp)
\\&=&
\frac{\partial^2}{\partial x^2}_{|x=1}\left(
x(1-x)^{2n}\right)\cdot\Tr (\square_q^{-s}P^\perp)
\\&=&
\left\{
\begin{array}{ll}
0&\mbox{if }n>1\\
2\sum_{\mu\in\Lambda^\vee\setminus\{0\}} \|\mu\|^{-2s}&\mbox{if }n=1\,\,.
\end{array}
\right.
\end{eqnarray*}

\end{document}